\documentclass[journal]{IEEEtran}
\usepackage{color} 
\usepackage{ifpdf}

\usepackage{cite}

\ifCLASSINFOpdf
  \usepackage[pdftex]{graphicx}

\else
 
\fi

\usepackage{graphicx}
\usepackage{epstopdf}
\usepackage{subfigure}
\usepackage{setspace}
\usepackage{amsmath}
\usepackage{amsfonts}
\usepackage{amssymb}
\ifCLASSOPTIONcompsoc
  \usepackage[caption=false,font=normalsize,labelfont=sf,textfont=sf]{subfig}
\else
  \usepackage[caption=false,font=footnotesize]{subfig}
\fi

\usepackage{fixltx2e}
\usepackage{footnote}
\makesavenoteenv{tabular}

\ifCLASSOPTIONcaptionsoff
  \usepackage[nomarkers]{endfloat}
 \let\MYoriglatexcaption\caption
 \renewcommand{\caption}[2][\relax]{\MYoriglatexcaption[#2]{#2}}
\fi

\usepackage{comment}
\usepackage{amsthm}
\usepackage{ifthen}

\newtheorem{definition}{Definition}

\newboolean{showcomments}
\setboolean{showcomments}{false}
\newcommand{\fliu}[1]{\ifthenelse{\boolean{showcomments}}
	{ \textcolor{red}{(FL:  #1)}}{}}
\newcommand{\gjp}[1]{\ifthenelse{\boolean{showcomments}}
	{ \textcolor{blue}{(Gjp:  #1)}}{}}

\def\ba{\begin{array}}
	\def\ea{\end{array}}
\newcommand{\beq}{\begin{equation}}
\newcommand{\eeq}{\end{equation}}
\newcommand{\bq}{\begin{eqnarray}}
\newcommand{\eq}{\end{eqnarray}}
\newcommand{\bqn}{\begin{eqnarray*}}
	\newcommand{\eqn}{\end{eqnarray*}}
\newcommand{\bee}{\begin{enumerate}}
	\newcommand{\eee}{\end{enumerate}}
\newcommand{\bi}{\begin{itemize}}
	\newcommand{\ei}{\end{itemize}}

\newcommand{\etab}{\boldsymbol{\eta}}

\begin{document}
	
	\title{Quantifying the Influence of Component Failure Probability on Cascading Blackout Risk}
	
	\author{Jinpeng~Guo,~\IEEEmembership{Student Member,~IEEE, }
		Feng~Liu,~\IEEEmembership{Member,~IEEE, }
        Jianhui~Wang,~\IEEEmembership{Senior Member,~IEEE, }
        Ming~Cao,~\IEEEmembership{Senior Member,~IEEE, }        
		Shengwei~ Mei,~\IEEEmembership{Fellow,~IEEE}
		\thanks{Manuscript received XXX, XXXX; revised XXX, XXX. \textit{(Corresponding author: Feng Liu)}.}
       	\thanks{Jinpeng Guo, Feng Liu, and Shengwei Mei are with the State Key Laboratory of Power Systems, Department of Electrical Engineering, Tsinghua University, Beijing, 100084, China  (e-mail: guojp15@mails.tsinghua.edu.cn, lfeng@tsinghua.edu.cn, meishengwei@tsinghua.edu.cn). } 
       	\thanks{Jianhui Wang is with the Department of Electrical Engineering at Southern Methodist University, Dallas, TX, USA and the Energy Systems Division at Argonne National Laboratory, Argonne, IL, USA (email: jianhui.wang@ieee.org)}
       	\thanks{Ming Cao is with the Institute of Engineering and Technology (ENTEG) at the University of Groningen, the Netherlands (email: m.cao@rug.nl)}	
        }
	
	\markboth{REPLACE THIS LINE WITH YOUR PAPER IDENTIFICATION NUMBER (DOUBLE-CLICK HERE TO EDIT}%
	{Shell \MakeLowercase{\textit{et al.}}:  Bare Demo of IEEEtran.cls for Journals}
	
	\maketitle
	
	\begin{abstract}
The risk of cascading blackouts greatly relies on failure probabilities of individual components in power grids. To quantify how component failure probabilities (CFP) influences blackout risk (BR), this paper proposes a sample-induced semi-analytic approach to characterize the relationship between CFP and BR. To this end, we first give  a generic component failure probability function (CoFPF) to describe CFP with varying parameters or forms. Then the exact relationship between BR and CoFPFs is built on the abstract Markov-sequence model of cascading outages. Leveraging a set of samples generated by blackout simulations, we further establish a sample-induced semi-analytic mapping between the unbiased estimation of BR and CoFPFs. Finally, we derive an efficient algorithm that can directly calculate the unbiased estimation of BR when the CoFPFs change. Since no additional simulations are required, the algorithm is computationally scalable and efficient. Numerical experiments well confirm the theory and the algorithm.

\end{abstract}
	
		
	\IEEEpeerreviewmaketitle

	\section{Introduction}
	\IEEEPARstart{D}{uring} the process of a cascading outage in power systems, the propagation  of component failures may cause serious consequences, even catastrophic blackouts \cite{r1,r2}. For the sake of effectively mitigating blackout risk, a naive way is to reduce the probability of blackouts, or more precisely, to reduce failure probabilities of system components by means of maintenance or so. Intuitively, it can be readily understood that component failure probabilities (CFPs) have great influence on blackout risk (BR). However, it is not clear how to efficiently quantify such influence, particularly in large-scale power grids. 
	
 	Generally, CFP largely depends on characteristics of system components as well as working conditions of those components. Since working conditions, e.g., system states, may change during a cascading outage, CFP varies accordingly.
 	Therefore, to quantitatively characterize  the influence of CFP on BR, two essential issues need to be addressed.	On the one hand, an appropriate probability function of component failure which can consider changing working conditions should be  well defined first to depict CFP . In this paper, such a function  is referred to as \emph{component failure probability function} (CoFPF). On the other hand, the quantitative relationship between  CoFPF and BR should also be explicitly established.
	
	For the first issue, a few CoFPFs  have been built in terms of specific scenarios. \cite{r3} proposes an end-of-life CoFPF of power transformers taking into account the effect of load conditions.  \cite{r4} considers the process and mechanism of tree-contact failures of transmission lines, based on which an analytic formulation is proposed. Another CoFPF of transmission lines  given in \cite{r5} adopts an exponential function to depict the relationship between CFP and some specific indices which can be calculated from monitoring data of transmission lines. This formulation can be extended to other system components in addition to transmission lines, e.g., transformers \cite{r6}. Similar works can be found in \cite{r7,r8}. In addition, some simpler CoFPFs are deployed in cascading outage simulations \cite{r9,r10,r11,r12}. Specifically, in the OPA model\cite{r9}, CoFPF is usually chosen as a monotonic function of the load ratio, while a  piecewise linear function is deployed in the hidden failure model\cite{r10}. Such simplified formulations have been widely used in various blackout models\cite{r11,r12}. However, it is worthy of noting that these CoFPFs are formulated for specific scenarios. In this paper, we adopt a generic CoFPF to facilitate establishing  the relationship between CFP and BR. 

    The second issue, i.e., the quantitative relationship between CoFPF and BR, is the main focus of this paper. Generally, since various kinds of uncertainties during the cascading process make the number of possible propagation paths explode with the increase of system scale, it is extremely difficult, if not impossible, to calculate the exact value of BR in practice. Therefore, analyticaly expressing the relationship between CoFPF and BR  is really challenging.  In this context, estimated BR by statistics, which is based on a set of samples generated by cascading simulations, appears to be the only practical substitute. 
    
    Among the sample-based approaches in BR estimation, Monte Carlo simulation (MCS) is the most popular one to date. However, MCS usually requires a  large number of samples with respect to specific CoFPFs for achieving satisfactory accuracy of estimation. Due to the intrinsic inefficiency,  MCS is greatly limited in practice, particularly in large-scale systems \cite{r13}. More importantly, it cannot  explicitly reveal the relationship between CoFPF and BR in an analytic manner. Hence, whenever parameters or forms of CoFPFs change, MCS must be completely re-conducted to generate new samples for correctly estimating BR, which is extremely time consuming. Moreover, due to the inherent strong nonlinearity between CoFPF and BR, when multiple CoFPFs change simultaneously, which is a common scenario in practice, BR cannot be directly estimated  by using the relationship between BR and individual CoFPFs. In this case, in order to correctly estimate BR and analyze the relationship, the required sample size will dramatically increase compared with the case that a single CoFPF changes. That indicates even efficient variance reduction techniques  (which may effectively reduce the sample size in a single scenario\cite{r14,r15}) are employed, the computational complexity will remain too high to be tractable. 
     
     The aforementioned issue gives rise to an interesting question: when one or multiple CoFPFs change, could it be possible to accurately estimate BR without re-conducting the extremely time-consuming blackout simulations? To answer this question, the paper proposes a sample-induced semi-analytic method to quantitatively characterize the relationship between CoFPFs and BR based on a given sample set. Main contributions of this work are threefold: 
    \begin{enumerate}
    \item Based on a generic form of CoFPFs, a cascading outage is formulated as a Markov sequence with appropriate transition probabilities. Then an exact relationship between BR and CoFPFs results.  
    \item Given a set of blackout simulation samples,  an unbiased estimation of BR is derived, rendering a semi-analytic expression of the mapping between BR and CoFPFs.   
    \item A high-efficiency algorithm is devised to directly compute BR when CoFPFs change, avoiding re-conducting any blackout simulations. 
    \end{enumerate} 
	
	The rest of this paper is organized as follows. In Section II, a generic formulation of CoFPF, an abstract model of cascading outages  as well as the exact relationship between CoFPF and BR are presented. Then the sample-induced mapping between the unbiased estimation of BR and CoFPFs is explicitly established in Section III. A high-efficiency algorithm  is presented in Section IV. In Section V, case studies are given. Finally, Section VI concludes the paper with remarks. 
  
\section{ Relationship between CoFPFs and BR }
The propagation of a cascading outage is a complicated dynamic process, during which many practical factors are involved, such as hidden failures of components,  actions of the dispatch/control center, etc. In this paper, we focus on the influence of random component failures (or more precisely, CFP) on the BR, where a cascading outage can be simplified into a sequence of component failures with corresponding system states, and usually emulated by steady-state models\cite{r9,r10,r12}. In this case, individual component failures are only related to the current system state while independent of previous states, which is known as the Markov property. This property enables an abstract model of cascading outages with a generic form of CoFPFs, as we explain below. 

\subsection{A Generic Formulation of CoFPFs}
To describe a CFP varying along with the propagation of cascading outages, a CoFPF is usually defined 
in terms of working conditions of the component. In the literature, CoFPF has various  forms with regard to specific scenarios \cite{r3,r4,r5,r6,r7,r8,r9,r10,r11,r12}. To generally depict the relationship between BR and CFP with varying parameters or forms, we first define an abstract CoFPF here. Specifically, the CoFPF of component $k$, denoted by $\varphi_k$, is defined as
\bq \label{eq1}
\varphi_k(s_k,\etab_k):=\mathbf{Pr}(\text{component}\; k\; \text{fails at}\; s_k \; \text{given} \; \etab_k) 
\eq
In \eqref{eq1}, $s_k$ represents the current working condition of component $k$, which can be load ratio, voltage magnitude, etc.  $\etab_k$ is the parameter vector. Both $\etab_k$ and the form of $\varphi_k$ represent the characteristics of the component $k$, e.g., the type and age of component $k$ \fliu{I don't understand this sentence}. It is worthy of noting that  the working condition of component $k$ varies during a  cascading outage, resulting in changes of the related CFPs. On the other hand, whereas the cascading process usually does not change $\etab_k$ and the form of $\varphi_k$, they can also be influenced due to controlled or uncontrolled factors,  such as maintenance and extreme weather, etc. In this sense, Eq. \eqref{eq1} provides a generic formulation to depict such properties of CoFPFs. 

\subsection{Formulation of Cascading Outages}	
		\begin{figure}[!t]
			\centering
			\includegraphics[width=0.95\columnwidth]{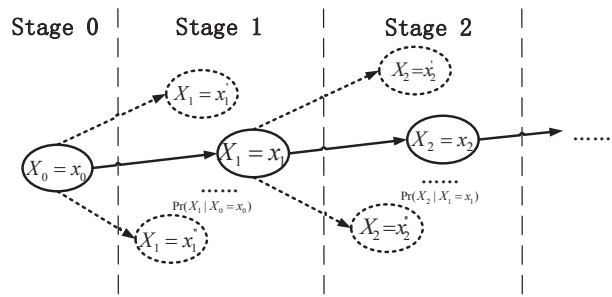}
			\caption{Cascading process in power systems}
			\label{fig.2}
		\end{figure}

\fliu{How about  we directly cite our previous paper and simplify the description?} 
In this paper, we are only interested in the paths of cascading outages that lead to blackouts as well as the associated load shedding. Then, according to \cite{r14}, cascading outages can be abstracted as a Markov sequence with appropriate transition probabilities.  
Specifically, denote $j\in \mathbb{N}$ as the sequence label and $X_j$ as the system state  at stage $j$ of a cascading outage. Here, $X_j$ can represent the power flow of transmission lines, the ON/OFF statuses of components, or  other system status of interest. The complete state space, denoted by $\mathcal{X}$, is spanned by all possible system states.  $X_0$ is the initial state of the system, which is assumed to be deterministic in our study. Under this condition, $\mathcal{X}$ is finite when only the randomness of component failures is considered \cite{r14}.  Note that $X_j$ specifies the working conditions of each component at the current stage (stage $j$) in the system, consequently determines $s_k$ in the CoFPF, $\varphi_k(s_k, \etab_k)$.
 Then an $n$-stage cascading outage can be represented by a series of states, $X_j$ (see Fig. \ref{fig.2})  and mathematically defined as below.

	\begin{definition} 
	\label{cascadingoutage}
	An $n$-stage cascading outage is a Markov sequence $
	Z:= \{ {X_0, X_1, ..., X_j, ..., X_n}$,  $X_j\in \mathcal{X}, \forall j\in \mathbb{N}  \}$
	with respect to a given joint probability series $g(Z)=g(X_{n} ,\cdots , X_{1} ,X_{0} )$.
	\end{definition}
 
In the definition above, $n$ is the total number of cascading stages, or the length of the cascading outage. Particularly, we denote the set of all possible paths of cascading outages in the power system by $\mathcal{Z}$. Since $\mathcal{X}$ and the total number of components are finite, $\mathcal{Z}$ is finite as well, albeit it may be huge in practice.\fliu{should we explain $\mathcal{Z}$ is finite?}\gjp{I prefer not mention this point in the original version} Then for a specific path, $z\in \mathcal{Z}$, of cascading outages, we have
\bqn
z&=&\{x_0, x_1, ..., x_j, ..., x_n\}\\
g(Z=z)&=&g(X_0=x_0, ..., X_j=x_j, ..., X_n=x_n)
\eqn
where $g(Z=z)$ is the joint probability of the path, $z$. For simplicity, we denote $g(Z=z)$ as $
g(z)=g({x_n}, \cdots, {x_1},{x_0})
$.

 Invoking the conditional probability formula and Markov Property, $g(z)$ can be further rewritten as 
\begin{equation} \label{eq8}
\begin{array}{rcl}
	g(z) & = &g({x_n}, \cdots, {x_1},{x_0}) \\
	&= &g_n({x_n}|{x_{n - 1}} \cdot  \cdot  \cdot {x_0})\cdot g_{n-1}({x_{n - 1}}|{x_{n - 2}} \cdots {x_0}) \\
	 & & \cdots g_1(x_1|x_0)\cdot g_0({x_0}) \\
	&= &g_n({x_n}|{x_{n - 1}})\cdot g_{n-1}({x_{n - 1}}|{x_{n - 2}}) \cdot  \cdot  \cdot g_0({x_0})
\end{array}
\end{equation} 
where
\begin{equation} \nonumber
	g_{j+1}({x_{j+1}}|{x_{j}})=\mathbf{Pr}(X_{j+1}=x_{j+1}|X_{j}=x_{j})
\end{equation} 
		\begin{equation} \nonumber
		g_0({x_0})=\mathbf{Pr}(X_0=x_0)=1
		\end{equation}

It is worthy of noting that this formulation is a mathematical abstraction of the cascading processes in practice and simulation models considering physical details \cite{r9,r10,r12}. Different from the high-level statistic models\cite{r18,r19}, it can provide an analytic way to depict the influence of many physical details, e.g., CFP, on the cascading outages and BR, which will be elaborately explained later on.  

\subsection{Formulation of Blackout Risk}
In the literature, blackout risk has various definitions \cite{r14,r20,r21} \fliu{refs?}. Here we adopt the  widely-used one, which is defined with respect to load shedding caused by cascading outages. 

Due to the intrinsic randomness of  cascading outages, the load shedding, denoted by $Y$, is also a random variable up to the path-dependent propagation of cascading outages. Therefore, $Y$ can be regarded as a function of cascading outage events, denoted by $Y:=h(Z)$.  Then the BR with respect to $g(Z)$ is defined as the expectation of the load shedding greater than a given level, $Y_0$. That is
		\begin{equation} \label{eq10}
		R_g(Y_0)=\mathbb{E}(Y\cdot\delta_{\{Y\ge{Y_0}\}})
		\end{equation}
where, $R_g(y_0)$ stands for the BR with respect to $g(Z)$ and $Y_0$; $\delta_{\{Y\ge{Y_0}\}}$ is the indicator function of $\{Y\ge{Y_0}\}$, given by
\begin{equation*} 
\delta_{\{Y\ge{Y_0}\}}:=\left\{\begin{array}{lll}
1 & &\text{if}\quad Y\ge{Y_0};\\
0 & &\text{otherwise}.
\end{array} \right.
\end{equation*}
In Eq. \eqref{eq10}, when the load shedding level is chosen as $Y_0=0$, it is simply the traditional definition of blackout risk. If $Y_0>0$, it stands for the risk of cascading outages with quite serious consequences, which is closely related to the renowned risk measures, value at risk(\textit{VaR}) and conditional value at risk(\textit{CVaR})\cite{r16}. Specifically, the risk defined in \eqref{eq10} is equivalent to \textit{CVaR}$_\alpha$ times $(1-\alpha)$ with respect to \textit{VaR}$_\alpha=Y_0$ with a confidence level   of $\alpha$.  
\subsection{Relationship Between BR and CoFPFs}
We first derive  the probability of cascading outages based on the generic form of CoFPFs. Then we characterize the relationship between BR and CoFPFs.

At stage $j$,  the working condition of component $k$ can be represented as a function of the system state $x_j$, denoted by $\phi_k(x_j)$. That is $s_k:=\phi_k(x_j)$. Hence the CFP of component $k$ at stage $j$ is $\varphi_k(\phi_k(x_j),\etab_k)$. For simplicity, we avoid the abuse of notion by letting $\varphi_k(x_j)$  to stand for $\varphi_k(\phi_k(x_j),{\etab_k})$. 

Considering  stages $j$  and $(j+1)$, we have

\begin{equation} \label{eq4}
g_{j+1}(x_{j+1}|x_{j})=\prod\limits_{k \in {F(x_{j})}} {{\varphi _k}( x_j  )} \cdot \prod\limits_{k \in {\bar{F}(x_{j})}} {(1 - {{\varphi _k}( x_j  )})}
\end{equation} 

In Eq. \eqref{eq4}, $F(x_{j})$ is the component set consisting of the components that are defective at $x_{j+1}$ but work normally at $x_{j}$, while $\bar{F}(x_{j})$ consists of components that work normally at $x_{j+1}$. \fliu{please check the previous sentence to make sur it is exact.} \gjp{$\bar{F}(x_{j})$ is not the complement of $F(x_{j})$, since there are components breaking down at former stages.}
With \eqref{eq4}, Eq. \eqref{eq8} can be rewritten as		
		\begin{equation} \label{eq2}
		\begin{array}{rcl}
		g(z)&=&{\prod\limits_{j=0}^{n-1}{g_{j+1}(x_{j+1}|x_{j})}}\\
		&=&\prod\limits_{j=0}^{n-1}\left[\prod\limits_{k \in {F( x_{j} )}} {{\varphi_{k}}({x_j})} \cdot \prod\limits_{k \in {\bar{F}( x_{j} )}} {(1 - {\varphi _{k}}({x_j}))} \right]\\		
		\end{array}
		\end{equation}		
	
Furthermore, substituting $Y=h(Z)$  into \eqref{eq10} yields
		\begin{equation} \label{eq11}
		\begin{array}{rcl}
		R_g(y_0) & = & \mathbb{E}(h(Z)\cdot\delta_{\{h(Z)\ge{Y_0}\}})\\
		&=&\sum\limits_{z \in \mathcal{Z}} {g(z)h(z){\delta _{\{ h(z) \ge {Y_0}\} }}} \\
		\end{array}
		\end{equation}

	Theoretically, the  relationship  between BR and CoFPFs can be established immediately by  	substituting \eqref{eq2} into \eqref{eq11}. However, this relationship cannot be directly applied  in practice, as we explain. Note that, according to \eqref{eq2}, the component failures occurring at different stages on a path of  cascading outages are correlated with one another. As a consequence, this long-range coupling, unfortunately,  produces  complicated and nonlinear correlation between BR and CoFPFs. In addition, since the number of components in a power system usually is quite large, the cardinality of $\mathcal{Z}$ can be  huge. Hence it is practically impossible to  accurately calculate BR with respect to the given CoFPFs by directly using \eqref{eq2} and \eqref{eq11}. To circumvent this problem, next we turn to using an unbiased estimation of BR as a surrogate, and propose a sample-based semi-analytic method to characterize the relationship.

\section{Sample-induced Semi-analytic Characterization}
\subsection{Unbiased  Estimation of BR}
To estimate BR, conducting MCS is the easiest and the most extensively-used way. The first step is to generate independent identically distributed (i.i.d.) samples of cascading outages and corresponding load shedding with respect to the joint probability series, $g(Z)$. Unfortunately,  $g(Z)$ is indeed unknown in practice. In such a situation, one can heuristically sample the failed components at each stage of possible cascading outage paths in terms of the corresponding system states and CoFPFs. Afterwards, system states at the next stage are determined with the updated system topology. This process repeats until there is no  new failure happening anymore. Then a path-dependent sample is generated. This method essentially  carries out sampling  sequentially using the \emph{conditional component probabilities} instead of the \emph{joint probabilities}. 
Eq. \eqref{eq2} provides this method with a mathematical interpretation, which is a  application of the Markov property of cascading outages. 


Suppose $N$ i.i.d. samples of cascading outage paths are obtained with respect to $g(Z)$. Let $Z_g:=\{z^i,i=1,\cdots,N\}$ record the set of these samples. Then, the $i$-th cascading outage path contained in the set is expressed by $z^i:=\{x_0^i,\cdots,x_{n^i}^i\}$, where $n^i$ is the number of  total stages of the $i$-th sample. For each $z^i \in Z_g$, the associated load shedding is given by $y^i=h(z^i)$. All $y^i$ make up the set of load shedding with respect to $g(Z)$, denoted by  $Y_g:=\{y^i,i=1,\cdots,N\}$. Then the unbiased estimation of BR is formulated as
\begin{equation} \label{eq12}
\hat{R}_g(Y_0)=\frac{1}{N}\sum\limits_{i = 1}^N {{y^i}{\delta _{\{ {y^i} \ge {Y_0}\} }}}=\frac{1}{N}\sum\limits_{i = 1}^N {{h(z^i)}{\delta _{\{ {h(z^i)} \ge {Y_0}\} }}}
		\end{equation}

Note that Eq. \eqref{eq12} applies to $g(Z)$ or the corresponding CoFPFs. That implies the underlying relationship between BR and CoFPFs relies on samples. Hence, whenever parameters or forms of the CoFPFs change, all samples need to be re-generated to estimate the BR, which is extremely time consuming,  even practically impossible. Next we derive a semi-analytic method by building a mapping between CoFPFs and the unbiased estimation of BR.  

\subsection{Sample-induced Semi-Analytic Characterization}
Suppose the samples are generated with respect to $g(Z)$. Then the sample set is $Z_g$,  and the set of load shedding is $Y_g$. When changing $g(Z)$ into $f(Z)$ (both are defined on $\mathcal{Z}$), usually all samples of cascading outage paths need to be regenerated. However, inspired by the sample treatment in Importance Sampling \cite{r14}, it is possible to avoid sample regeneration by revealing the underlying relationship between  $g(Z)$ and $f(Z)$. 
Specifically, for a given  path $z^i$, we define
\begin{equation}
\label{eq:w_def}
w(z^i):=\frac{f(z^i)}{g(z^i)}, \quad (z^i\in\mathcal{Z})
\end{equation} 
\fliu{should we use $z^i$?} then each sample in terms of $f(z)$ can be represented as the sample of $g(z)$ weighted by $w(z)$. Consequently, the unbiased estimation of BR in terms of $f(Z)$ can be directly obtained from the sample generated with respect to $g(Z)$, as we explain. 

From Eqs. \eqref{eq12} and \eqref{eq:w_def}, we have 
\begin{equation} \label{eq13}
\hat{R}_f(Y_0)=\frac{1}{N}\sum\limits_{i = 1}^N {w(z^i){h(z^i)}{\delta _{\{ {h(z^i)} \ge {Y_0}\} }}}
\end{equation}
Obviously, when $w(z)\equiv 1,z\in\mathcal{Z}$, \eqref{eq13} is equivalent to \eqref{eq12}. Moreover, Eq. \eqref{eq13} is an unbiased estimation  by noting that
\begin{equation} \label{eq14}
\begin{array}{rcl}
\mathbb{E}(\hat{R}_f(Y_0) )& = & \mathbb{E}\left(\frac{f(Z)}{g(Z)}{h(Z)}{\delta _{\{ {h(Z)} \ge {Y_0}\} }}\right) \\
& = & \sum\limits_{z\in\mathcal{Z}}{ g(z) \times \frac{f(z)}{g(z)}{h(z)}{\delta _{\{ {h(z)} \ge {Y_0}\} }}  }\\
& = & R_f(Y_0)
\end{array}
\end{equation}
Note that in Eqs. \eqref{eq13} and \eqref{eq14}, only the information of $h(z)$ is required. One crucial implication is that the BR with respect to $f(Z)$, $R_f$, can be estimated directly with no need of  regenerating cascading outage samples. This feature can further  lead to an efficient algorithm to analyze BR under varying CoFPFs, which  will be discussed next.

\section{Estimating BR with Varying CoFPFs}
\subsection{Changing a Single CoFPF}
We first consider a simple case, where a single CoFPF changes. Suppose CoFPF of the $m$-th component changes from $\varphi_m$ to $\bar{\varphi}_m$\footnote{For simplicity, we use the notion $\bar{\varphi}_m$ to denote the new CoFPF of component $m$, which may have a new function form or parameters $\etab_m$. }, and the corresponding joint probability series changes from $g(Z)$ to $f(Z)$. Considering a sample of cascading outage path  generated with respect to $g(Z)$, i.e., $z^i \in Z_g$, we have
\begin{small}		
	\begin{equation} \label{eq15}
	\begin{array}{rcl}
	f(z^i)&=&{\prod\limits_{j=0}^{n^i-1}{f_{j+1}(x_{j+1}^i|x_{j}^i)}}\\
	&=&{\prod\limits_{j=0}^{n^i-1}{\left[\prod\limits_{k \in {F^m( x_{j}^i) }} {{\varphi _k}({x_j^i})}\cdot \prod\limits_{k \in {\bar{F}^m( x_{j}^i )}} {(1 - {\varphi _{k}}({x_j^i}))}\right]}  }\cdot \\
	&&{\cdots \Gamma(\bar{\varphi}_{m},z^i) }
	\end{array}
	\end{equation}
\end{small}
where,
\begin{small}	
	\begin{equation} \label{eq7}
	\Gamma(\bar{\varphi}_{m},z^i) = \left\{ {\begin{array}{*{20}{c}}
		\prod\limits_{j=0}^{n_m^i-1}{(1 - {\bar{\varphi} _m}({x_j^i}))}&:& \text{if}\quad n_m^i=n^i\\
		{\bar{\varphi} _m}({x_{n_m^i}^i}) \prod\limits_{j=0}^{n_m^i-1}{(1 - {\bar{\varphi} _m}({x_j^i}))}&:& \text{otherwise}\\
		\end{array}} \right.
	\end{equation}
\end{small}
Here, $n_m^i$ is the stage at which the $m$-th component experience an outage. Particularly, $n_m^i=n^i$ when the $m$-th component is still working normally  at the last stage of the cascading outage path. Component set $F^m( x_{j}^i ):=F( x_{j}^i )\setminus\{m\}$ consists of all the elements in $F( x_{j}^i )$ except for $m$. Similarly $\bar{F}^m( x_{j}^i ):=\bar{F}( x_{j}^i )\setminus\{m\}$ is  the component set including  all the elements in $\bar{F}( x_{j}^i )$ except for $m$ \fliu{we should explain what is $F_m(x^i_j)$ and $\bar{F}_m(x^i_j)$}.   According to \eqref{eq:w_def}, the sample weight is 
\begin{equation} \label{eq16}
w(z^i)=\frac{f(z^i)}{g(z^i)}=\frac{\Gamma(\bar{\varphi}_{m},z^i)}{\Gamma({\varphi}_{m},z^i)}
\end{equation}

Substituting \eqref{eq16} into \eqref{eq13}, the unbiased estimation of BR is 
\begin{equation} \label{eq17}
\hat{R}_f(Y_0)=\frac{1}{N}\sum\limits_{i = 1}^N {\frac{\Gamma(\bar{\varphi}_{m},z^i)}{\Gamma({\varphi}_{m},z^i)}{h(z^i)}{\delta _{\{ {h(z^i)} \ge {Y_0}\} }}}	
\end{equation}

Eq. \eqref{eq17} provides an unbiased estimation of BR after changing a CoFPF by only using the original samples. 
\subsection{Changing Multiple CoFPFs}
In this section, we consider the general case that multiple CoFPFs change simultaneously. Invoking the expression of $f(z^i)$ in \eqref{eq15}, we have  

\begin{equation} \label{eq19}
g(z^i)=\prod\limits_{k\in K}{\Gamma(\varphi_k,z^i)}
\end{equation} 
\begin{equation} \label{eq20}
f(z^i)=\prod\limits_{k\in K_c}{\Gamma(\bar{\varphi}_k,z^i)} \cdot \prod\limits_{k\in K_u}{\Gamma(\varphi_k,z^i)}
\end{equation}
where, $K$ is the complete set of all  components in the system; $K_c$ is the set of components whose CoFPFs change;  $K_u$ is the set of others,  i.e.,   $K=K_c\cup K_u$; $\bar{\varphi}_k$ is the new CoFPF of the $k$-th component. 		

According to \eqref{eq19} and \eqref{eq20}, the sample weight is given by     
\begin{equation} \label{eq21}
w(z^i)=\prod\limits_{k\in K_c}{ \frac{ {\Gamma(\bar{\varphi}_k,z^i)} }{\Gamma({\varphi}_k,z^i)}}			
\end{equation}
Substituting \eqref{eq21} into \eqref{eq13} yields 
\begin{equation} \label{eq22}
\hat{R}_f(Y_0)=\frac{1}{N}\sum\limits_{i = 1}^N { \left( \prod\limits_{k\in K_c}{ \frac{ {\Gamma(\bar{\varphi}_k,z^i)} }{\Gamma({\varphi}_k,z^i)}}  {h(z^i)}{\delta _{\{ {h(z^i)} \ge {Y_0}\} }} \right) }			
\end{equation}  

Eq. \eqref{eq22} is a generalization  of Eq. \eqref{eq17}.  \eqref{eq22} provides a mapping between the unbiased estimation of BR and CoFPFs. When multiple CoFPFs change, the unbiased estimation of BR can be directly calculated by using \eqref{eq22}. Since no additional cascading outage simulations are required,  and only algebraic calculations are involved, it is  computationally  efficient.      

\subsection{Algorithm} 
To clearly illustrate the algorithm, we first rewrite \eqref{eq22} in a matrix form as
\begin{equation} \label{eq23}
\hat{R}_f(Y_0)=\frac{1}{N}\mathbf{L}\mathbf{F_p}			
\end{equation}   
In \eqref{eq23}, $\mathbf{L}$ is an $N$-dimensional row vector, where $\mathbf{L}_i={h(z^i)}{\delta _{\{ {h(z^i)} \ge {Y_0}\} }}/g(z^i)$. $\mathbf{F_p}$ is an $N$-dimensional column vector, where $\mathbf{F_p}_{i}=f(z^i)$. 

We further define two $N \times k_{a}$ matrices $\mathbf{A}$ and $\mathbf{B}$, where $k_{a}$ is the total number of all components, $\mathbf{A}_{ik}=\Gamma(\varphi_k,z^i)$, $\mathbf{B}_{ik}=\Gamma(\bar{\varphi}_k,z^i)$. According to \eqref{eq20},  we have
\begin{equation} \label{eq24}
\mathbf{F_p}_i = \prod\limits_{k\in K_c}{\mathbf{B}_{ik}} \cdot \prod\limits_{k\in K_u}{\mathbf{A}_{ik}}
\end{equation}

Then the algorithm  is given as follows.

\noindent\rule[0.25\baselineskip]{0.5\textwidth}{1pt}	
\begin{itemize}
	\small
	\item {\bf Step 1: Generating samples.} Based on the system and blackout model in consideration, generate a set of i.i.d. samples. Record the sample sets $Z_g$ and  $Y_g$, as well as the row vector $\mathbf{L}$.    
	\item {\bf Step 2:  Calculating $\mathbf{F_p}$.} Define the new CoFPFs for each component in $K_c$, and calculate $\mathbf{B}$ and $\mathbf{A}$. Then calculate $\mathbf{F_p}$ according to \eqref{eq24}. Particularly, instead of calculation, $\mathbf{A}$ can be saved in Step 1.
	\item {\bf Step 3:  Data analysis.} According to  \eqref{eq23}, estimate the blackout risk for the changed CoFPFs.
\end{itemize}
\noindent\rule[0.25\baselineskip]{0.5\textwidth}{1pt}


\subsection{Implications}
The proposed method has important implications in blackout-related analyses. Two of typical examples are:  efficient estimation of blackout risk considering extreme weather conditions and the risk-based maintenance scheduling.

For the first case, as well known,  extreme weather conditions (e.g., typhoon) often occur for a short time but affect a wide range of components. The failure probabilities of related components may increase remarkably. In this case, the proposed method can be applied to fast evaluate the consequent risk in terms of the weather forecast.

For the second case, since maintenance  can considerably reduce CoFPFs, the  proposed method allows an efficient identification of the  most effective  candidate  devices in the system for mitigating blackout  risks.  Specifically, suppose that one only considers simultaneous maintenance of  at most $d_m$ components. Then the number of possible scenarios is up to $\sum\nolimits_{d=1}^{d_m}{C({k_{a}},d)}$, which turns to be intractable in a large practical system ($C({k_{a}},d)$ is the number of $d$-combinations from $k_a$ elements). Moreover, in each scenario, a great number of cascading outage simulations are required to estimate the blackout risk, which is extremely time consuming, or even practically impossible. In contrast, with the proposed method, one only needs to generate the sample set in the base scenario. Then the blackout risks  for  other  scenarios  can be  directly calculated  using only algebraic calculations, which is very simple and computationally    efficient. 

\section{Case Studies}
\subsection{Settings}
In this section, the numerical experiments are carried out on the IEEE 300-bus system with a simplified OPA model omitting the slow dynamic \cite{r9,r17}. Its basic sampling steps are summarized as follows.
\newline
\noindent\rule[0.25\baselineskip]{0.5\textwidth}{1pt}	
\begin{itemize}
	\small
	\item {\bf Step 1: Data initialization.} Initialize the system data and parameters. Particularly, define specific CoFPFs for each component. The initial state is $x_0$.
	\item {\bf Step 2: Sampling outages.} At stage $j$ of the $i$-th sampling, according to the system state, $ x_j^i $ , and CoFPFs, simulate the component failures with respect to the failure probabilities. 
	\item {\bf Step 3:  Termination judgment.} If new failures happen in Step 2, recalculate the system state $ x_{j+1}^i $ at stage $j+1$ with the new topology, and go back to Step 2; otherwise, one sampling ends. The corresponding samples are $z^i = \{x_0, x_1^i \cdots  x_j^i\}  $ and $y^i$. If all $N$ simulations are completed, the sampling process ends. 
\end{itemize}
\noindent\rule[0.25\baselineskip]{0.5\textwidth}{1pt}

In this simulation model, the state variables $X_j$ are chosen as the ON/OFF statuses of all components and power flow on corresponding components at stage $j$. Meanwhile, the random failures of transmission lines and power transformers are considered. The total number of them is $k_{a}=411$, and the CoFPF we use is 

\begin{equation} \label{eq26}
\varphi_k (s_k,\etab_k)= \left\{ {\begin{array}{*{20}{c}}
	{p_{\min }^k}&:&{s^k < s_d^k}\\
	{\frac{{p_{\max }^k - p_{\min }^k}}{{s_u^k - s_d^k}}(s^k - s_d^k) + p_{\min }^k}&:&{else}\\
	{p_{\max }^k}&:&{s^k > s_u^k}
	\end{array}} \right.
\end{equation}
where $s^k$ is the load ratio of component $k$ and $\etab_k=[p^k_{min},p^k_{max},s^k_d,s^k_u]$. Specifically,  $p^k_{min }$ denotes the minimum failure probability of component $k$ when the load ratio is less than $s^k_d$; $p^k_{\max } $ denotes the maximum failure probability when the load ratio is larger than $s^k_{u} $. Usually it holds $0<p^k_{min}<p^k_{max}<1$, which  depicts certain probability of hidden failures in protection devices. Particularly, we use the following initial parameters to generate $Z_g$: $s_d^k=0.97$, $s_u^k=1.3$, $p^k_{max }=0.9995$, $p^k_{min } \sim U[0.002,0.006]$, $k\in K$. 

It is worthy of noting that the simulation process with specific settings mentioned above is a simple way to emulate the propagation of cascading outages. We only use it to demonstrate the proposed method. The proposed method  can apply when  other more realistic models and parameters are adopted.   

\subsection{Unbiasedness of the  Estimation of  Blackout Risk}
\begin{figure}[!t]
	\centering
	\includegraphics[width=0.95\columnwidth]{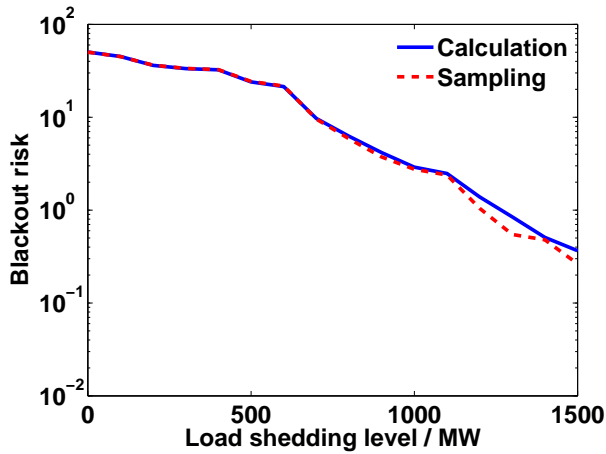}
	\caption{BR Estimation  with sampling and calculation}
	\label{fig.1}
\end{figure}

In this case, we will show that our method can  achieve unbiased  estimation  of blackout risk. We first carry out 100,000 cascading outage simulations with the initial parameters. Then the sample set $Z_g$ and related $\mathbf{L}$ are obtained. Afterward we randomly choose a set of failure components, $K_c$, including two components. Accordingly, we modify the  parameters, $\etab_k$, of their CoFPFs to $\bar{p}^k_{min }=p^k_{min}-0.001,k \in K_c$. In terms  of the new settings and various load shedding levels, we estimate the blackout risks by using \eqref{eq23}. For comparison, we re-generate  100,000 samples under the new settings and estimate the BRs by using \eqref{eq12}. The results are given in Fig. \ref{fig.1}.

Fig. \ref{fig.1} shows that the blackout risk estimations with two methods are almost the same, which indicates the proposed method can achieve unbiased estimation of BR. Note that  our  method requires no more simulations, which is much more efficient than traditional MCS, and scalable  for  large-scale systems. 

\subsection{Parameter Changes in CoFPFs}
In this case, we test the performance of our method when parameters, $\etab_k$, of some CoFPFs change. The sample set $Z_g$ and related $\mathbf{L}$ are based on the 100,000 samples with respect to the initial parameters. We consider two different settings: 1) $Y_0=0$; 2) $Y_0=1500$. Statistically, we  obtain  $\hat{R}_g(0)=54.2$ and $\hat{R}_g(1500)=0.43$, respectively. 

Since there are 411 components in the  system, we consider  $k_a=411$ different scenarios. Specifically, we change each CoFPF individually  by letting  $\bar{p}^k_{min }=p^k_{min}-0.001$. Using the  proposed method, the blackout risks can be  estimated quickly. Some results  are presented  in  Table \ref{t1} ($Y_0=0$MW) and Table \ref{t2} ($Y_0=1500$MW). Particularly, the average computational times of unbiased estimation of BR in each scenario in Table \ref{t1} and \ref{t2} are 0.004s and 0.00001s, respectively.\fliu{we should give the computational time for each case in Tables I and II}  

\begin{table}[htp]\footnotesize
	\caption{BR estimation when parameters of CoFPFs change  $(Y_0=0)$ }
	\label{t1}
	\centering
	\begin{tabular}{c|ccc}
		\hline \hline
		Component & $p^k_{min}$ & Blackout  &Risk reduction \\
		index&$(\times 10^{-3})$&risk& ratio$(\%)$ \\
		\hline
		204 &5.8&50.97&6.0\\
		312 &5.1&51.93&4.2\\
		372 &3.0&53.43&1.5\\
		114 &4.1&53.46&1.4\\
		307 &3.7&53.53&1.3\\
		410 &3.6&53.68&1.0\\
		117 &4.0&53.69&1.0\\
		63 &2.3&53.70&1.0\\
		259 &3.0&53.90&0.6\\
		126 &2.7&53.92&0.5\\
		...&...&...&...\\
		Mean value &4.0&54.13&0.2\\
		\hline \hline
	\end{tabular}
\end{table}

\begin{table}[htp]\footnotesize
	\caption{BR estimation when parameters of CoFPF changes $(Y_0=1500)$ }
	\label{t2}
	\centering
	\begin{tabular}{c|ccc}
		\hline \hline
		Component & $p^k_{min}$ & Blackout  &Risk reduction \\
		index&$(\times 10^{-3})$&risk&ratio$(\%)$\\
		\hline
		259 &3.0&0.358&16.7\\
		254 &4.3&0.380&11.4\\
		204 &5.8&0.382&11.0\\
		312 &5.1&0.403&6.0\\
		93  &3.0&0.403&6.0\\
		325 &2.2&0.404&5.9\\
		48  &3.0&0.404&5.9\\
		378 &2.6&0.405&5.7\\
		116 &2.2&0.406&5.5\\
		305 &2.2&0.407&5.3\\
		...&...&...&...\\
		Mean value &4.0 &0.427&0.5\\
		\hline \hline
	\end{tabular}
\end{table}

Table \ref{t1} shows the top ten scenarios having the  lowest value of BR as well as the average risk of all scenarios. Whereas reducing   failure probabilities of certain components   can effectively mitigate the blackout  risk, others have little  impact. For example, decreasing  $p^k_{min}$ of component $\#$204 results in $6.0\%$ reduction of blackout risk, while the average ratio of risk reduction is only $0.2\%$. This result implies that   there may exist some critical components which play a core role in  the propagation of cascading outages and promoting load shedding. Our method enables a scalable way to efficiently identify those components, which may facilitate  effective  mitigation of blackout risk.    

\begin{figure}[htp]
	\centering
	\includegraphics[width=1\columnwidth]{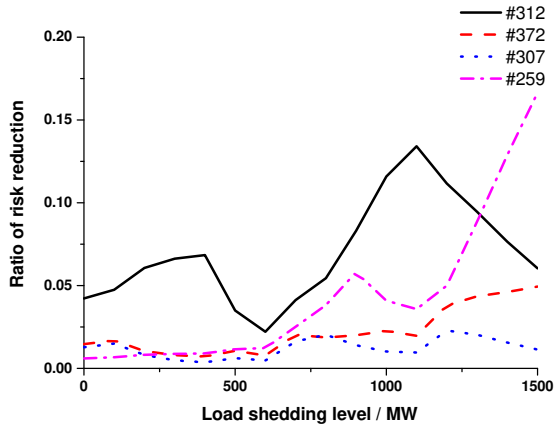}
	\caption{Ratio of risk reduction under different load shedding levels}
	\label{fig.3}
\end{figure}

When considering $Y_0=1500$MW, which is associated with serious blackout events, it is interesting to  see  the most  influential components in Table \ref{t2} are different from that in Table \ref{t1}. In other words, the impact of CoFPF varies with different load shedding levels, which demonstrates the complex relationship between BR and CoFPFs.  To better show this point, we  choose four typical components ($\#$312, $\#$372, $\#$307 and $\#$259) and calculate the risk reduction ratios with respect to a series of load shedding levels. As shown in Fig. \ref{fig.3}, component $\#259$ has little influence on medium to small size of blackouts, while considerably reducing the risk of large-size blackouts. It implies this component has a significant contribution to the propagation of cascading outages.  
 In contrast, some other components, such as $\#307$ and $\#372$, result in similar risk reduction ratios for various load shedding levels. They are likely to cause certain direct load shedding, but have little influence on the propagation of cascading outages. In terms of component  $\#312$, the curve of risk reduction ratio exhibits a multimodal feature, which means 
the changes of such a CoFPF may have much larger influence on BR for some load shedding levels than others.
All these results demonstrate a highly complicated  relationship between BR and CoFPFs. Our method indeed provides a  computationally efficient tool to conveniently analyze such relationships in practice. 

\begin{table}[htp]\footnotesize
	\caption{Risk reduction ratio when parameters of CoFPFs change $(\%)$ }
	\label{t3}
	\centering
	\begin{tabular}{cc|cc}
		\hline \hline
		Component&Risk reduction ratio&Component&Risk reduction ratio\\
		index&$(Y_0=0)$&index&$(Y_0=1500)$\\
		
		\hline
		(204,312) &10.2&(204,259) &25.5\\
		(204,372) &7.5 &(254,259) &24.5\\
		(114,204) &7.4 &(93,259)  &22.6\\
		(204,307) &7.3 &(259,312) &22.2\\
		(63,204)  &7.0 &(259,378) &21.8\\
		(204,410) &7.0 &(259,305) &21.8\\
		(117,204) &7.0 &(259,325) &21.6\\
		(204,259) &6.6 &(116,259) &21.3\\
		(126,204) &6.6 &(40,259)  &21.1\\
		(204,301) &6.5 &(48,259)  &21.1\\
		...       &... & ...      &... \\
		Mean value&0.3 &Mean value&1.0 \\
		\hline \hline
	\end{tabular}
\end{table}
Then we  decrease $p^k_{min}$ of two CoFPFs simultaneously.  $Z_g$, $\mathbf{L}$  are the same as before, and the number of scenarios is $C(k_m,2)=C(411,2)$. The calculated unbiased estimations with respect to two $Y_0$ are shown in Table \ref{t3}.
It is no surprise that the risk reduction ratios are more remarkable compared with the  results in Table \ref{t1} and  Table \ref{t2}, where only one CoFPF decrease. However, it should be noted that the pairs of components in  Table \ref{t3} are not simple combinations of the ones shown in Table \ref{t1}/Table \ref{t2}. The reason is that the relationship between BR and CFP is complicated and nonlinear, which actually results in the difficulties in  analyses as we demonstrate in Section III.

\subsection{Form Changes in CoFPFs }
In this case, we show the performance of the proposed method when the form of CoFPFs changes. The new form of CoFPF of component $k$ is modified into
\begin{equation} \label{eq27}
\bar{\varphi}_k (s_k,\etab_k) = \left\{ {\begin{array}{*{20}{c}}
	{max(\varphi_k,ae^{bs^k})}&:&{s^k < s_u^k}\\
	{p_{\max }^k}&:&{s^k \ge s_u^k}
	\end{array}} \right.
\end{equation}
where $a=p_{\min }^k$ and $b=\frac{ln(p_{\max }^k/p_{\min }^k)}{s_u^k}$ are corresponding parameters. $\varphi_k$, $p_{\min }$, $p_{\max }$, $s_u^k$, $s_d^k$ are the same as the ones in \eqref{eq26}. Obviously, the failure probabilities  of  individual components  are amplified in this case.

\begin{figure}[!t]
	\centering
	\includegraphics[width=0.95\columnwidth]{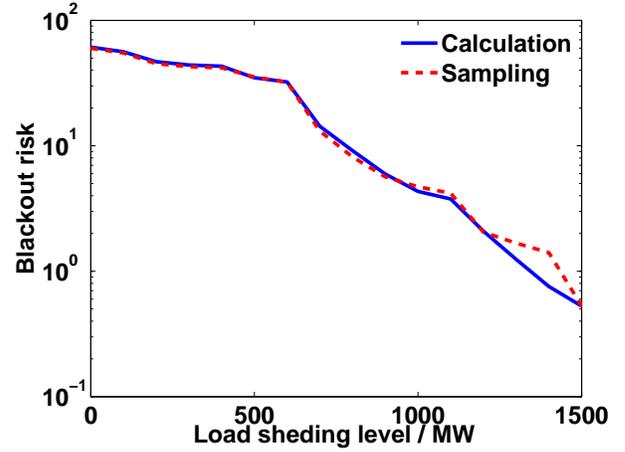}
	\caption{Estimation of the blackout risk with sampling and calculation}
	\label{fig.4}
\end{figure}

Similar to the previous cases, we use the proposed method to calculate the unbiased estimations of blackout risks  when some CoFPFs change from \eqref{eq26} to \eqref{eq27}. The comparison results of our method and MCS are presented in Fig. \ref{fig.4}. In addition, the blackout risk in several typical scenarios is shown in Tabs \ref{t5} and \ref{t7}. These results empirically confirm the efficacy of our method.  

\begin{table}[!t]\footnotesize
	\caption{Risk increase ratio when form of CoFPF changes $(\%)$ }
	\label{t5}
	\centering
	\begin{tabular}{cc|cc}
		\hline \hline
		Component&Risk increase ratio&Component&Risk increase ratio\\
		index&$(Y_0=0)$&index&$(Y_0=1500)$\\
		\hline
		204 &12.3& 204 &22.1\\
		312 &5.5 & 259 &21.0\\
		259 &0.8 & 254 &16.2\\
		264 &0.7 & 312 &7.5\\
		372 &0.6 & 264 &6.8\\
		...       &... & ...      &... \\
		Mean value&0.1 &Mean value&0.2 \\
		\hline \hline
	\end{tabular}
\end{table}

\begin{table}[!t]\footnotesize
	\caption{Risk increase ratio when forms of CoFPFs change $(\%)$ }
	\label{t7}
	\centering
	\begin{tabular}{cc|cc}
		\hline \hline
		Component&Risk increase ratio&Component&Risk increase ratio\\
		index&$(Y_0=0)$&index&$(Y_0=1500)$\\
		\hline
		(204,312) &17.8 &(204,259) &48.7\\
		(204,259) &13.1 &(254,259) &43.5\\
		(204,264) &12.9 &(204,254) &42.3\\
		(204,372) &12.9 &(204,312) &30.7\\
		(204,307) &12.8 &(201,204) &30.2\\
		...       &... & ...      &... \\
		Mean value&0.2 &Mean value&0.5 \\
		\hline \hline
	\end{tabular}
\end{table}

\subsection{Computational Efficiency }
We carry out all tests on a computer with an Intel Xeon E5-2670 of 2.6GHz and 64GB memory. It takes 107 minutes to generate $Z_g$ (100,000 samples) with respect to the initial parameters. Then we enumerate all $\sum\nolimits_{d=1}^{2}{C({k_{a}},d)=84666}$ scenarios. In each scenario, the parameters and forms of one or two CoFPFs are changed (cases in Part C and D, Section V, respectively).  Then we calculate the unbiased estimations of BR with $Y_0=0$ and $Y_0=1500$ by using the proposed method  \fliu{I don't understand what do  ``Case C'' and ``Case D'' in Table V stand for...}.   The complete computational times are given in Tables \ref{t6} and \ref{t8}. It is worthy of noting that the computational times of risk estimation in tables are the total times for enumerating \emph{all} the 84666 scenarios. It takes only about 10 minutes to calculate $\mathbf{B}$, and additional 10 minutes to computing BRs. On the contrary, it will take  about $107\times 84666 \approx 9,000,000$ minutes for the MCS method,  which is practically intractable. 

\begin{table}[!t]\footnotesize
	\caption{Computing time when parameter or form changes in CoFPF (Min.)}
	\label{t6}
	\centering
	\begin{tabular}{c|ccc}
		\hline \hline
		    &  Calculate $\mathbf{B}$  & Estimate risk $(Y_0=0)$ &  Total\\
		\hline
		Parameter change &10.0&9.2 &19.2\\
		Form change &10.7&9.3 &20.0\\
		\hline \hline
	\end{tabular}
\end{table}
 
 \begin{table}[!t]\footnotesize
	\caption{Computing time when parameter or form changes in CoFPF (Min.)}
	\label{t8}
	\centering
	
		\begin{tabular}{c|ccc}
		\hline \hline	
		    & Calculate $\mathbf{B}$  & Estimate risk $(Y_0=1500)$ &  Total\\
		\hline
		Parameter change &10.0&0.01 &10.01\\
		Form change &10.7&0.01 &10.71\\
		\hline \hline
	\end{tabular}
	
\end{table}
\section{Conclusion with Remarks}
In this paper, we propose a sample-induced semi-analytic method to efficiently quantify the influence of  CFP on BR. Theoretical analyses and numerical experiments show that:
\begin{enumerate}
\item With appropriate formulations of cascading outages and a generic form of CoFPFs, the relationship between CoFPF and BR is exactly characterized and can be effectively estimated based on samples.

\item Taking advantages of the semi-analytic expression between CoFPFs and the unbiased estimation of BR, the BR can be efficiently estimated when CoFPFs change, and the relationship between the CFP and BR can be analyzed correspondingly.
\end{enumerate}

Numerical experiments reveal that the long-range correlation among component failures during cascading outages is really complicated, which is often considered closely related to self-organized criticality and power low distribution. Our results serve as a step towards providing a scalable and efficient tool for further understanding  the failure correlation and the mechanism of cascading outages.  Our ongoing work  includes the application of the proposed method in making maintenance plans and  risk evaluation considering extreme weather conditions, etc. 
	
\fliu{The style of reference is not correct. We need to check it carefully. }

\end{document}